\documentstyle[12pt]{article}
\begin{document}

\begin{flushright}
RIBS-PH-4/97   \\

q-alg/9703032
\end{flushright}

\begin{center}
{\Large \bf  Green Function on the q-Symmetric Space $SU_q(2)/U(1)$}
\end{center}

\vspace{1cm}

\begin{flushleft}
H. Ahmedov$^1$  and I. H. Duru$^{1,2}$

\vspace{.5cm}

{\small
1. TUBITAK -Marmara Research Centre, Research Institute for Basic
Sciences, Department of Physics, P.O. Box 21, 41470 Gebze, 
Turkey\footnote{E--mail: 
ahmedov@mam.gov.tr and duru@mam.gov.tr.}.

2. Trakya University, Mathematics Department, P.O. Box 126, 
Edirne, Turkey.}
\end{flushleft}

\vspace{2cm}

\begin{center}
{\bf \small  Abstract}
\end{center} 

Following the introduction of the invariant distance on the non-commutative 
C-algebra of the quantum group $SU_q(2)$, the Green function and 
the Kernel on the q-homogeneous space $M=SU(2)_q/U(1)$  are derived. 
A path integration is formulated. Green function for the free massive
scalar field on the non-commutative 
Einstein space $R^1\times M$ is presented.

\vspace{2cm}
\begin{center}
March 1997
\end{center}

\pagebreak

\vfill
\eject

\section{Introduction}
How the quantum mechanical effects should be altered  if we replace the space-time  
continuum with a non-commutative geometry is an exciting question. To  answer this
question we have to formulate the known quantum mechanical problems over the 
non-commutative spaces. Since we lack satisfactory mathematical tools, construction  
 of the Schr\"{o}dinger equations over the non-commutative spaces  is difficult 
 and sometimes arbitrary \cite{kn:1}. First of all when we do not have a differentiable 
 manifold it is problematic to find the correct operators replacing the derivatives. 
 If however the non-commutative geometry is given as a quantum group space this problem may be 
 solved in a natural way for having the action of the q-algebra generators 
 one is not required to deal with q-differential calculus. If this is not the case,
 since it is always possible to build up an integration theory on a given set,
 the path integrals may in principle be the suitable method of 
 quantization for the non-commutative geometries in general.
 Therefore the derivation of the Green functions over the  non-commutative 
 spaces is of interest. 
 We should also remember that in the usual quantum physics defined over the 
 commutative spaces the Green functions which are the vacuum (or temperature) 
 expectation values of two point field operators play important role \cite{kn:2}.
 Even for the q-group spaces the construction of the Green 
 functions seems to be important step in the formulation of many q-deformed 
 quantum mechanical problems. 

 The experience we have in the derivation of Green functions over the (undeformed)
 group manifolds is quite rich; and it is also well known that many 
 non-relativistic quantum mechanical problems  are related to 
 the particle motion over these manifolds \cite{kn:3}. 
 Therefore we hope that constructing  the Green functions on the q-group spaces 
 may lead  meaningful definitions of these problems over the non-commutative 
 geometries.It is important to stress that if we know the formulations of the 
 non-relativistic potential problems  we also gain insight in some field 
 theoretical effects. In fact, the calculations of many field theoretical 
 problems like the pair creations in the given cosmologies or in the external
 electromagnetic fields, and the Casimir interactions may formally become 
 equivalent to some non-relativistic potential problems. For example to 
 investigate the pair production in the Robertson-Walker space-time expanding with the 
 factor $a(t)$ one has to calculate the Green function of a particle moving 
 in the one-dimensional potential $V(t)=a^{-2}(t)$ with the time $t$ playing 
 the  role of coordinate \cite{kn:4}.

 Motivated by the considerations summarized above the construction of
 the Green functions  over the quantum group spaces is the subject of this work.
 The specific example we study is the quantum symmetric space 
 $M=SU_q(2)/U(1)$.Quantum symmetric spaces have already been subject of some 
 interesting investigations.
 For example the well known relations between the special functions and the classical
 groups have been generalized to the quantum groups \cite{kn:5}; and 
 quantum spheres are studied \cite{kn:6}. Recently the homogeneous space of 
 $E_q(2)$ is considered and q-Schr\"{o}dinger equation on it is constructed 
 \cite{kn:7}.

 In section II after a brief review of the invariant distance concept, 
 we outline a method for  constructing the Green functions by two classical 
 group examples.In this method, which is applicable to the quantum groups,
 one first constructs the one-point ``Green function", then obtains the 
 Green function depending on two points by the group action.
 
 In Section III the invariant distance for the quantum group $A=Pol(SU_q(2))$
 and the q-symmetric space $M=SU_q(2)/U(1)$ is constructed and its 
 properties are demonstrated.

 In Section IV the one-point ``Green function" is derived on the space $M$,
from which we obtain the Green function in section V.

 In Section VI we introduce the time development Kernel on the space $M$.
 Having this Kernel in hand the non-commutative path integration is also formulated. 

Finally in Section VII the Green function for the massive scalar field on 
$R^1\times M$ which is the non-commutative version of the 
Einstein space  is constructed.
 
The basic definitions and the established results about the Hopf algebra A
which we use in our work are given in the Appendix.

\vfill
\eject

\renewcommand{\theequation}{II.\arabic{equation}}
\setcounter{equation}{0}

\section{ Method for Constructing Green Functions. 
Examples from the classical Lie groups }
The Green function of the free particle motion over a Lie group 
manifold and its homogeneous spaces depend on the invariant distance 
between two points. When we attempt to construct Green functions over the 
quantum group spaces the first problem we have to face is the introduction  
of the invariant distance. To overcome this problem it is instructive to 
review the case of the classical Lie groups for the purpose of developing 
a method for constructing the Green functions which can also be 
employed for the quantum groups.We briefly study two  examples.

We first consider the real line $\cal R$ which is 
the homogeneous space with respect to the  translation group $T(x)$. 
The Green function of the free particle motion    
over $\cal R$  depends on the invariant distance $\mid x-x'\mid$ and 
on the momentum $ p\in(-\infty,\infty)$ which is the weight of the unitary 
irreducible representation.  
The homogeneity of $\cal R$ under the action of $T(x)$ implies
\begin{equation}
 {\cal G}^p(x,x') =T^{-1}(x'){\cal G}^p(x,0).
\end{equation}
The above equation suggests  that the invariant  Green function   on the homogeneous space
can be obtained by the group action if the one-point ``Green function" 
${\cal G}^p(x,0)$ is known.

As the second example we consider the classical symmetric space $SU(2)/U(1)$.
If an element $g\in SU(2)$ is parametrized as $g=kak'$ with 
\begin{equation}
k=  \left(
\begin{array}{cc}
e^{i\psi /2}  &  0  \\
0  &  e^{-i\psi /2}
\end{array}
\right ),   
a=  \left(
\begin{array}{cc}
\cos\theta /2  & i\sin\theta /2   \\
i\sin\theta /2   &  \cos\theta /2
\end{array}
\right )   
\end{equation}
the symmetric space M which is topologically equivalent to $S^2$ is 
represented as
\begin{equation}
M=g\sigma(g^{-1}).
\end{equation}
Here $\sigma$ is the involutive automorphism having the property
\begin{equation}
\sigma (k)=k, \ \ \ \ \ \sigma (a)=a^{-1}.
\end{equation}
The Green function  ${\cal G}$ on M depends on the group invariants which 
are the weight of the irreducible representation $l$=0,1,2,... and on 
the invariant distance given by
\begin{equation}
\rho (x_1,x_2)=1-\frac{1}{2}Tr(x_1x_2^{-1}) 
\end{equation}
with $x_1,x_2\in M$.
Using the group property $g_2^{-1}g_1=g_{12}$ we can write
\begin{eqnarray}
Tr(x_1x_2^{-1})= Tr(g_1\sigma(g_1^{-1})(g_2\sigma(g_2^{-1}))^{-1})=
Tr(g^{-1}_2g_1\sigma((g^{-1}_2g_1)^{-1}))=  \nonumber \\
=Tr(g_{12}\sigma(g_{12}^{-1}))=Tr(x_{12}).
\end{eqnarray}
If we fix one of the points as 
$x_1=\left (  \begin {array} {cc}
1 & 0 \\
0 & 1
\end {array} \right )$, we obtain a formula depending on only one point.
Taking the advantage of this formula we can first construct a one-point 
``Green function", then by the action of the group element we arrive at the 
Green function which is dependent on two points. The equation satisfied 
by the one-point ``Green function" is
\begin{equation}
({\cal C}-(l+1/2)^2){\cal G}^l(x)=\delta (x)
\end{equation}
where ${\cal C}$ is the center of the enveloping algebra $U(su(2))$.
Once we obtain the solution of the above equation we can derive the 
Green function simply by the group action as
\begin{equation}
{\cal G}^l(x_1,x_2)=T(g_2^{-1}){\cal G}^l(x_1)={\cal G}^l(g_2x_1\sigma(g_2^{-1}))
\end{equation}
which satisfies
\begin{equation}
({\cal C}-(l+1/2)^2){\cal G}^l(x_1,x_2)=\delta (x_1-x_2).
\end{equation}

\renewcommand{\theequation}{III.\arabic{equation}}
\setcounter{equation}{0}

\section{An Invariant Distance on the Quantum Group $SU_q(2)$ }

The coordinate functions $a,a^{*},b,b^{*}$ of the Hopf algebra $A$ 
(see App.A) satisfy the commutation relations \cite{kn:6}:
\begin{eqnarray}
b^{*}b=bb^{*}, \  \ ba=qab, \ \ b^{*}a=qab^{*},  \nonumber  \\
aa^{*}+b^{*}b=1, \ \ a^{*}a+q^2bb^{*}=1
\end{eqnarray}
The $*$-representation of the quantum group $A$ (the representation of the 
C-algebra) in the Hilbert space ${\cal H}$ with the orthonormal basis 
$\{ \mid n \rangle \}$, n=0,1,2,... with $q<1$ is given by \cite{kn:6}
\begin{eqnarray}
a\mid n\rangle=(1-q^{2n})^{1/2}\mid n-1\rangle, \ \ \ \ \ \
b\mid n\rangle=e^{i\phi}q^n\mid n\rangle, \nonumber \\
b^{*}\mid n\rangle=e^{-i\phi}q^{n+1}\mid n\rangle, \ \ \ \ \ \
a^{*}\mid n\rangle=(1-q^{2n+2})^{1/2}\mid n+1\rangle
\end{eqnarray}
is irreducible for fixed  $\phi\in(0,2\pi)$.

In  $q\rightarrow 1$ limit the relations (III.2) define the three dimensional 
space which is the manifold of $SU(2)$.To understand the situation 
in the deformed case we recall the usual quantum mechanics in which the
physical systems are defined by the vectors in the Hilbert space
and the self adjoint operators correspond to  the observables including  
coordinates. In the same manner for quantum groups we construct 
 self-adjoint operators $X$ from the linear combinations 
of the coordinate functions given in (III.2) and define the expectation values
of these operators as the points of the quantum group space as 
\begin{equation}
X_\psi=\langle\psi\mid X\mid \psi \rangle, \ \ \ X\in A, \psi\in{\cal H}. 
\end{equation}
The above definition can be carried to co-product space $X\otimes Y$ with
$X,Y\in A$ to define the invariant distance on $SU_q(2)$ which should go
to the corresponding classical limit as $q\rightarrow 1$ and should have 
the following properties:
\begin{equation}
\rho (g_1,g_2)=\rho (g_2,g_1),
\end{equation}
\begin{equation}
\rho (g_1g,g_2g)=\rho (g_2,g_1),
\end{equation}
\begin{equation}
\rho (g_1,g_2)>0,
\end{equation}
\begin{equation}
\rho (g,g)=0,
\end{equation}
\begin{equation}
\rho (g_1,g_2)<\rho (g_1,g_3)+\rho(g_3,g_2).
\end{equation}

Motivated by formula  (II.5) of the previous section we suggest the 
following hermitian operator in $A\otimes A$ :
\begin{equation}
R=(\tau\otimes S) \Delta
(1-\frac{1}{[2]_q}Tr_q(d^{1/2}))
\end{equation}
Here $d^{1/2}$ is
the matrix of the unitary irreducible corepresentation of the Hopf algebra 
$A$  with  weight $1/2$ (see App.A). $[\cdot ]_q$ is define as
$[x]_q=\frac{q^x-q^{-x}}{q-q^{-1}}$ and the q-trace is
given by
\begin{equation}
Tr_q(d^{1/2})=\sum_{j=-1/2}^{1/2}q^{-2j}(d^{1/2}_{jj}).
\end{equation} 
 $\tau$ is the automorphism of $A$ defined as
\begin {equation}
\tau(d^{1/2})
=\left(
\begin{array}{cc}
q^{-1}a  &  b  \\
-qb^{*}  &  qa^{*}
\end{array}
\right ).
\end{equation}  
Expectation value of the operator $R$ of (III.9) in the Hilbert space 
$H\otimes H$
defines the correct invariant distance for A. In fact 
this operator possesses  all the properties of (III.4--III.8). 

(i) Symmetry condition of (III.4) is fulfilled for
\begin{equation}
\sigma R=R.
\end{equation}
Here $\sigma$ is the flip homomorphism 
$\sigma (x\otimes y)=y\otimes x$; $x,y\in A$.

(ii) Invariance condition of (III.5) takes the form of
\begin{equation}
\langle\Delta\otimes \Delta R \rangle_4=R
\end{equation}
where  $\langle\cdot\rangle_4$  is the map from 
$A\otimes A\otimes A\otimes A$ into $A\otimes A$ is given by
\begin{equation}
\langle a\otimes b\otimes c\otimes d\rangle_4=
(a\otimes c)\phi(b\tau^{-1}(d))
\end{equation}
with $\phi(x)$ being the invariant integral on the quantum group $A$
(see app.B).

(iii) The operator $R$ is positive. 
To show this we have to construct the 
basis in the Hilbert space $H\otimes H$ in which it is diagonal. 
We choose the eigenfunctions  as
\begin{equation}
\psi_l=\sum_{j=-l}^lv^l_j\mid l+j\rangle\otimes \mid l-j\rangle .
\end{equation}
We then sandwich the operator $R$ of (III.9)between the above states 
to get the recurrence relations for the unknown coefficients $v^l_j$
\begin{equation}
[2]_qv^l_j(1-\frac{2q^{2l+1}}{[2_q]}-E_l)=
v^l_{j-1}[l+j][l-j+1]+v^l_{j+1}[l-j][l+j+1].
\end{equation}
These coefficients are normalized as
\begin{equation}
\sum_{j=-l}^l\overline{v^l_j}v^l_j =1.  
\end{equation}
$E_l$ in (III.16) is the spectrum of the self-adjoint operator $R$ and 
$[x]$  is defined as $[x]=\sqrt{1-q^{2x}}$. 
As an example consider the eigenvalue and eigenfunction 
for  $l=0$ state 
\begin{equation}
\psi_0=\mid 0 \rangle \otimes \mid 0\rangle.  
\end{equation}
The corresponding eigenvalue 
\begin{equation}
E_0=\frac{q^{-1}-q}{q^{-1}+q}
\end{equation}
is positive for $q<1$. Similar demonstrations can be done for all other 
values of $l$ to prove that  $E_l$ is positive. 

(iv) condition of  (III.7) is also satisfied for
\begin{equation}
m(\tau^{-1}\otimes id) R=0
\end{equation}
where m is the operation of multiplication in the C-algebra and $\tau^{-1}$
is the inverse of the involution (III.11).

(v) Finally the triangular inequality of (III.8) reads
\begin{equation}
\langle\psi\mid (id\otimes \sigma)(R\otimes 1)\mid \psi\rangle <
\langle\psi\mid (R\otimes 1 +1\otimes R)\mid \psi\rangle , \ \ \ 
\mid \psi\rangle \in H\otimes H \otimes H.
\end{equation}
As in the case (iii) we can show that the self-adjoint operator $\Pi$ 
which is defined as
\begin{equation}
\Pi=R\otimes 1 +1\otimes R - (id\otimes \sigma)(R\otimes 1) 
\end{equation}
is positive. For example for the eigenfunction 
\begin{equation}
\psi_{0}=\mid 0\rangle\otimes\mid 0\rangle\otimes\mid 0\rangle, \ \ \
\langle \psi_0\mid\psi_0\rangle=1
\end{equation}
the corresponding  eigenvalue which is given by
\begin{equation}
E_0=\frac{q^{-1}-q}{q^{-1}+q}
\end{equation}
is positive for $q<1$.

Before closing this section we give the invariant distance on the coset  
space $M=A/K$ where $K=Pol(U(1))$(see App. B). 
We first introduce the involutive automorphism
\begin {equation}
\beta(d^{1/2})
=\left(
\begin{array}{cc}
a  &  -b  \\
-c  &  d
\end{array}
\right )
\end{equation}  
It is clear that it does not change the quantum subgroup $K$ (see App. A).
By virtue of this automorphism we introduce the coordinate functions on 
the q-symmetric space $M=A/K$ as
\begin {equation}
t_{ij}=(d^{1/2}\beta(S(d^{1/2})))_{ij}.
\end{equation}  
The operator of the invariant distance on $M$ is then given by
\begin{equation}
R_M=(\tau\otimes S)\Delta
(1-\frac{1}{[2]_q}Tr_q(t)).
\end{equation}
Introducing (III.26) into the above formula we obtain
\begin{equation}
R_M=(\tau\otimes S)\Delta\xi
\end{equation}
where $\xi=-q^{-1}bc $ is the element of the  two-sided coset space
$H=K/A/K$ (see App.A).

\renewcommand{\theequation}{IV.\arabic{equation}}
\setcounter{equation}{0}

\section{ One-Point ``Green Function" }
We will follow the method we introduced in Sec.II in the construction
of the Green function over $M$. We first have to obtain the  one-point 
``Green function" ${\cal G}_q^l(\xi)$ .It is defined by the 
 deformation of (II.7) as
\begin{equation}
({\cal C}-[l+\frac{1}{2}]^2) {\cal G}_q^l(\xi) =\delta_q(\xi)
\end{equation}
where ${\cal C}$ is the center of the Hopf algebra $U(su_q(2))$ (see App.C).
The invariant q-delta function $\delta_q (\xi)$ is a linear functional  over 
the two-sided coset  space $H$ which for any function $f\in A[0,0]$ 
satisfies
\begin{equation}
\langle \delta_q(\xi)\mid f(\xi)\rangle = f(0), \ \ \  f\in A[0,0]
\end{equation}
where  scalar product is the one given in App.B.
It is easy to verify that the q-delta function can be represented as
\begin{equation}
\delta_q(\xi)= \sum_{l=0}^{\infty} [l+\frac{1}{2}] d_{0,0}^l(\xi)
\end{equation}
where $d_{0,0}^l(\xi)$ is the q-zonal function \cite{kn:6}
\begin{equation}
d_{0,0}^l(\xi)= \phi_{2,1}(q^{-2l},q^{2(l+1)}, q^{2l}
\mid q^2,q\xi).
\end{equation}
We can also verify by direct substitution that  
the one-point ``Green function" of (IV.1) can be represented as
\begin{equation}
{\cal G}_q^l(\xi)=\sum_{n=0}^{\infty}  [n+\frac{1}{2}]
\frac{ D_{0,0}^l(\xi)}
{[n+\frac{1}{2}]^2-[l+\frac{1}{2}]^2}.
\end{equation}
The above summation can be executed to give an expression in terms of the  
q-hypergeometric function:
\begin{equation}
{\cal G}_q^l(\xi)=\gamma^l\xi^{-l-1}
\phi_{2,1}(q^{-2(l+2)},q^{-2(l+2)}, q^{-4(l+2)}
\mid q^{-2},q^{-2}\xi^{-1}) 
\end{equation}
where $\gamma^l$ is the normalization constant
\begin{eqnarray}
\gamma^l=\frac{1}{[l][l+1]} (\int_0^1 d_q \xi {\cal G}_q^l(\xi))^{-1}= 
\frac{\sqrt{q}}{(-1-l)_{-q}}  \nonumber \\
(\phi_{2,1}(q^{-2(l+2)},q^{-2l)}, q^{-4(l+2)}
\mid q^{-2},q^{-2}))^{-1} 
\end{eqnarray}
with $(x)_q=\frac{1-q^x}{1-q}$.

\renewcommand{\theequation}{V.\arabic{equation}}
\setcounter{equation}{0}

\section{Green Function over $M=A/K$ }
Having the one-point Green function in hand we can now introduce  
the  Green function ${\cal G}_q^l(M\otimes M)$ on the q-symmetric space 
$M=A/K$ as
\begin{equation}
 {\cal G}_q^l(M\otimes M)=(\tau\otimes S)\Delta{\cal G}_q^l(\xi).
\end{equation}
The equation satisfied by the above Green function is 
\begin{equation}
\left ( id\otimes \{ {\cal C} - [l+\frac{1}{2}]^2 \} \right )
\circ {\cal G}_q^l(M\otimes M)=
\delta_q(M\otimes M)
\end{equation}
where the invariant q-delta function  is given by
\begin{equation}
\delta_q(M\otimes M)=( \tau\otimes S ) \Delta \delta_q(\xi) .
\end{equation}
Substituting (IV.3) into the above equation and (IV.5) into (V.2) 
we obtain the following representations for the 
invariant  q-delta function and the Green function: 
\begin{equation}
\delta_q(M\otimes M)= \sum_{l=0}^{\infty}\sum_{j=-l}^l [l+\frac{1}{2}]
\tau(d_{0,j}^l(M))\otimes\overline{d_{0,j}^l}(G).
\end{equation}
\begin{equation}
{\cal G}_q^l(M\otimes M)= \sum_{n=0}^{\infty} \sum_{j=-n}^n[n+\frac{1}{2}]
\frac{\tau(d_{j,0}^n(M))\otimes \overline{d_{0,j}^n}(M)}
{[l+\frac{1}{2}]^2 - [n+\frac{1}{2}]^2}
\end{equation}
Using the representation (IV.6) of the one-point ``Green function" we have 
another expression for   ${\cal G}_q^l(M\otimes M)$ in terms of the q-hypergeometric 
function as
\begin{equation}
{\cal G}_q^l(M\otimes M)= \gamma^l (\tau\otimes S) \Delta
 (\xi^{-l-1}\phi_{2,1}(q^{-2(l+2)},q^{-2(l+2)}, q^{-4(l+2)}
\mid q^{-2},q^{-2}\xi^{-1}) )
\end{equation}
For any operator function $f(M)\in M$ and the linear operator $P$ of 
the dual Hopf algebra $U(su_q(2))$ (see App.C) the q-delta  function of (V.4)
satisfies 
\begin{equation} 
\langle \delta_q(M\otimes M) \mid id\otimes f\rangle_2 = f
\end{equation}
and
\begin{equation} 
\langle (P\otimes id) \delta_q(M\otimes M) \mid id\otimes f\rangle_2 = 
\langle\delta_q(M\otimes M) \mid id\otimes P^{*}f\rangle_2 = 
{\cal F}
\end{equation}
Here the inner product $\langle\cdot \mid\cdot\rangle_2 $
which is defined as
\begin{equation} 
\langle x_1\otimes x_2 \mid y_1\otimes y_2 \rangle_2= 
x_1 y_1\phi (x_2\tau^{-1}(y_2)); \ \ \ x_1,x_2,y_1,y_2\in M 
\end{equation}
is a map
\begin{equation} 
(M\otimes M)\times  (M\otimes M) \rightarrow M.
\end{equation}
Before closing this section we like to consider the inhomogeneous 
equation for a given constant $E$ and operator function $f(M)\in M$
\begin{equation}
( {\cal C} - E ) \circ {\cal F} = f.
\end{equation}
As in the classical case  the solution is obtained by using 
the Green function
\begin{equation}
{\cal F}  = F_{0}+
\langle {\cal G}_q^E(G\otimes G;E) \mid  f\otimes id \rangle_2
\end{equation}
where $F_{0}$ is the complete solution of the homogeneous equation.

\renewcommand{\theequation}{VI.\arabic{equation}}
\setcounter{equation}{0}

\section{ Kernel on the q-symmetric space $M$}
We introduce the unitary operator in terms of the real time interval 
$t$ and the center of the enveloping algebra ${\cal C}$ as 
\begin{equation}
U(t)=e^{it{\cal C}}
\end{equation}
satisfying the semigroup property
\begin{equation}
U(t)U(t')=U(t+t').
\end{equation}
The one-point ``q-Kernel" is then given by
\begin{equation}
{\cal K}_q(\xi ,t) = U(t) \delta_q(\xi )
\end{equation}
Inserting the  representation of the one-point q-delta function from (IV.3) 
we obtain
\begin{equation}
{\cal K}_q(\xi,t)= \sum_{l=0}^{\infty} [l+\frac{1}{2}]
e^{it[l+\frac{1}{2}]^2}d_{0,0}^l(\xi)
\end{equation}
which is connected to  the  one-point ``Green  function" through 
the relation
\begin{equation}
{\cal K}_q(\xi ,t) =\int_{-\infty}^{\infty} dE
e^{itE} {\cal G}_q^E(\xi).
\end{equation}
It is obvious that the above Kernel satisfies the equation
\begin{equation}
(i\partial_t + {\cal C}){\cal K}_q(\xi ,t) = 0.
\end{equation}

The two-point q-Kernel  is defined in a manner parallel to the definition 
of the one-point ``q-Kernel" as
\begin{equation}
{\cal K}_q(M\otimes M ,t) = (U(t)\otimes id)\delta_q(M\otimes M )
\end{equation}
Inserting the representation of the two-point q-delta function of (V.4)
into the above equation we have
\begin{equation}
{\cal K}_q(M\otimes M ,t)= \sum_{l=0}^{\infty}\sum_{j=-l}^l 
e^{it[l+\frac{1}{2}]^2}
[l+\frac{1}{2}]
\tau(d_{0,j}^l(M))\otimes \overline{d_{0,j}^l(M)}.
\end{equation}
The triple invariant product  $\langle\cdot\mid\cdot\rangle_3 $ defined by 
\begin{equation}
\langle x_1\otimes x_2\otimes x_3\mid y_1\otimes y_2\otimes y_3\rangle_3=
x_1y_1\otimes x_3y_3 \phi(x_2\tau^{-1}(y_2))  
\end{equation}
which is the  map 
\begin{equation}
(M\otimes M\otimes M)\times (M\otimes M\otimes M) \rightarrow 
M\otimes M
\end{equation}
enables us the derivation of the important  property of $K_q(M\otimes M,t)$:
\begin{equation}
\langle{\cal K}_q(M\otimes M ,t)\otimes 1 \mid  
1\otimes {\cal K}_q(M\otimes M ,t') \rangle_3= 
{\cal K}_q(M\otimes M ,t+t'). 
\end{equation}
Using this property  we can introduce the path integral representation 
for the Kernel on the non-commutative space. Indeed from (VI.11) follows 
\begin{equation}
{\cal K}_q(M\otimes M ,T)=
\langle{\cal K}_q(M\otimes M ,T/2)\otimes 1 \mid  
1\otimes {\cal K}_q(M\otimes M ,T/2) \rangle_3
\end{equation}
and
\begin{equation}
{\cal K}_q(M\otimes M ,T/2)=
\langle{\cal K}_q(M\otimes M ,T/4)\otimes 1 \mid  
1\otimes {\cal K}_q(M\otimes M ,T/4) \rangle_3
\end{equation}
Inserting (III.13) into (III.12) and using the short hand notation 
${\cal K}_q(T)= {\cal K}_q(M\otimes M ,T)$ we get 
\begin{equation}
{\cal K}_q(T)=
\langle\langle {\cal K}_q(T/4)\otimes 1 \mid 1\otimes {\cal K}_q(T/4) \rangle_3
 \otimes 1  \mid
\langle{\cal K}_q(T/4)\otimes 1 \mid 1\otimes {\cal K}_q(T/4) \rangle_3
\otimes 1 \rangle_3
\end{equation}
Continuing  the above process n times and  taking the limit 
$n\rightarrow \infty$ we arrive at a path integral formula:
\begin{equation}
{\cal K}_q(M\otimes M ,T)=\lim_{n\rightarrow \infty}
\{ ({\cal K}_q(M\otimes M ,T/2^n)\}_n 
\end{equation}
where $\{\cdot \}_n$ stands for $n$ times repeated $\langle\cdot\rangle_3$
map.

\renewcommand{\theequation}{VII.\arabic{equation}}
\setcounter{equation}{0}

\section{Green function for the massive free scalar field on the 
q-Einstein space $R^1\times M$}
For the commutative translation group parametrized by $t$, 
the commutative and co-commutative Hopf algebra of $Fun(t)$ is given by
\begin{equation}
\delta t=1\otimes t + t\otimes 1 , \ \ \ S(t)=-t, \ \ \ \varepsilon(t)=0.  
\end{equation}
 The one-point Kernel for the free particle motion on $(t,s)$ ``space-time" 
 is the usual one
 \begin{equation}
{\cal K}(t,s)=(-4i\pi s)^{-1/2}e^{-t^2/4s} .
\end{equation}
The Kernel on $(R^1\times M,s)$ is expressed as 
\begin{equation}
{\cal K}(\xi,t,s)={\cal K}(t,s){\cal K}_q(\xi,s) 
\end{equation}
where ${\cal K}_q(\xi,s)$ is given by (VI.4). Using the Schwinger-DeWitt 
representation 
\begin{equation}
{\cal G}(\xi,t;m^2)=-i\theta (t)\int_0^\infty e^{-im^2s}{\cal K}(\xi,t,s)ds
\end{equation}
with $Im (m^2)<0$ and $\theta (t)$ being the step function we get the 
one-point ``Green function" over $R^1\times M$ for the scalar field with 
the mass $m$. 
Performing the integration over $ds$ we obtain
\begin{equation}
{\cal G}(\xi,t;m^2)=\theta(t) \sum_{n=0}^\infty e^{i t\sqrt{[n+1/2]_q+m^2}}
\frac{[n+1/2]_q}{\sqrt{[n+1/2]_q^2+m^2}}d_{00}^n(\xi).
\end{equation}
The above Green function  satisfies
\begin{equation}
(\partial^2_t -{\cal C} +m^2){\cal G}(\xi,t;m^2)=\delta(t)\delta_q(\xi ).
\end{equation}
Following the procedure of  Sec. 5 we obtain the invariant Green function on the
space $R^1\times M$ depending on two points 
\begin{equation}
{\cal G}(Y\otimes Y,;m^2)=
(\tau\otimes S)\Delta {\cal G}(\xi,t;m^2) ,
\end{equation}
where $Y=t\times M$ and the operations $\Delta$ and $\tau$ on $t$ are given
by $\Delta t=\delta t$ and $\tau(t)=t$. 
The Green function (VII.7) satisfies 
\begin{equation}
(id\otimes (\partial^2_t -{\cal C} +m^2)){\cal G}(Y\otimes Y;m^2)=
\delta(t\otimes1-1\otimes t)\delta_q(M\otimes M )
\end{equation}

\vfill
\eject

\begin{center}
{\Large \bf Appendix }
\end{center}

\renewcommand{\theequation}{A.\arabic{equation}}
\setcounter{equation}{0}

{\large A. Hopf Algebra $A=Pol(SU(2)_q)$}

The algebra  of polynomials $A=A(SU_q(2))$  form the $\ast$Hopf algebra 
or real quantum group. The coordinate functions  $\pi_{ij}$ are given by
\begin {equation}
\pi_{ij}d^{1/2}=\pi_{ij}\left(
\begin{array}{cc}
a  &  b  \\
c  &  d
\end{array}
\right )=d^{1/2}_{ij} 
\end{equation}  
where $d^{1/2}$ is the matrix of the unitary irreducible corepresentation  
of the Hopf algebra $A$.
The coproduct $\Delta$, counit $\varepsilon$ and antipode $S$ acts as
\begin {equation}
\Delta \circ d^{1/2}_{ij}=d^{1/2}_{ik}\otimes d^{1/2}_{kj}
\end{equation}  
\begin {equation}
S(d^{1/2})=\left(
\begin{array}{cc}
d  &  -qb  \\
-q^{-1}c  &  d
\end{array}
\right ) 
\end{equation}  
\begin {equation}
\varepsilon\circ d^{1/2}=\left(
\begin{array}{cc}
1  &  0  \\
0  &  1
\end{array}
\right )
\end{equation}  
The $\ast$-operation  is
\begin {equation}
\left(
\begin{array}{cc}
a^{*}  &  b^{*}  \\
c^{*}  &  d^{*}
\end{array}
\right ) =
\left(
\begin{array}{cc}
d  &  -q^{-1}c  \\
-q b  &  a
\end{array}
\right ) 
\end{equation}  

\vspace{3mm}

\renewcommand{\theequation}{B.\arabic{equation}}
\setcounter{equation}{0}

{\large B. Harmonic Analysis on the  Coset Space $M=A/K$}

\vspace{2mm}

i. {\small The Cartan Decomposition of A.}

The quantum group $K=Pol(U(1))$ is the Hopf algebra with coordinate 
functions $t, \ t^{-1}$
\begin{equation}
\Delta_{U}\circ t^{\pm}=t^{\pm}\otimes t^{\pm}, \ \ S_{U}\circ t^{\pm}=t^{\mp}, \ \ 
\varepsilon_{U}\circ t^{\pm }=1 
\end{equation}

The Hopf algebra $K$ is the subalgebra of the quantum group A defined by the 
homomorphism 
\begin{equation}
\psi_k\left(
\begin{array}{cc}
a  &  b  \\
c  &  d
\end{array}
\right )=
\left(
\begin{array}{cc}
t  &  0  \\
0  &  t^{-1}
\end{array}
\right )
\end{equation}
The left and right unitary co-representation of $K$ in $A$ is given 
by homomorphisms 
\begin{equation}
L_u = (\psi_k\otimes id)\circ\Delta , \ \ \ 
R_u = (id\otimes \psi_k)\circ\Delta.
\end{equation}
The subspaces $A[j,i]$; $j,i\in Z$  defined by
\begin{equation}
A[j,i]=\{ x\in A : L_u \circ x = t^j\otimes x ; 
R_u \circ x = x \otimes t^i \}
\end{equation}
form the basis of the Hopf algebra A
\begin{equation}
A = \sum_{j,i\in  Z} \oplus A[j,i].
\end{equation}
The quantum coset space $M=A/K$ and   two-sided coset space $H=K/A/K$ are
the subspaces of $A$ 
\begin{equation}
M = \sum_{j,\in  Z} \oplus A[0,j], \ \ \ \ \ \ H=A[0,0]. 
\end{equation}

\vspace{3mm}

ii.{\small Harmonic analysis on $M$}

\vspace{2mm}

The irreducible co-representation of A is constructed in the space $V$ 
of homogeneous polynomials of degree $l$. 
\begin{equation}
T: \Delta V=V\otimes A
\end{equation}
The basis in $V$ is generated by the elements
\begin{equation}
e^l_j=a^{l+j}b^{l-j}.
\end{equation}
The matrix elements of the irreducible co-representation of
the $\ast$Hopf algebra are given in terms of $\xi=-q^{-1}bc\in H$ by
\begin{equation}
D^l_{0j}(M) =\lambda^l_j
\phi_{21}(q^{2(j-l)},q^{2(j+l+1)}q^{2(j+l)}\mid q^2,q\xi)c^jd^j; 
\ \ j=0,1,...,l 
\end{equation}
and
\begin{equation}
D^l_{0j}(M) =\lambda^l_j a^{-j}b^{-j}
\phi_{21}(q^{2(-j-l)},q^{2(-j+l+1)}q^{2(-j+l)}\mid q^2,q\xi); 
\ \ j=-l,-l+1,...,0 
\end{equation}
Here $\phi_{21}$ is the q-hypergeometric function 
\cite{kn:6}, and $\lambda^l_j$ is defined as
\begin{equation}
\lambda^l_j =q^{\mid j\mid (\mid j\mid -l)}
\left[
\begin {array} {c}
l \\
\mid j \mid
\end{array}
\right ]_{q^2}^{1/2}
\left[
\begin {array} {c}
l-\mid j\mid \\
\mid j \mid
\end{array}
\right ]_{q^2}^{1/2}
\end{equation}
The co-representation (B.7) is unitary with respect to the scalar product
\begin{equation}
\langle x \mid y\rangle =\psi(x^\ast y), \ \ \ \  x,y\in V
\end{equation}
where $\psi$ is the invariant integral  on $A$
\begin{equation}
\Psi(z)=\int_0^1 d\xi_q {\cal P}(z), \ \ \ \ \ \ \  z\in A
\end{equation}
and ${\cal P}$ is the projection operator ${\cal P}A[i,j]\rightarrow A[0,0]$.
With respect to the invariant integral (B.13) the matrix elements $D^l_{0j}$ 
satisfy the orthogonality condition 
\begin{equation}
\langle D^l_{0i}(M)\mid D^k_{0j}(M)\rangle=[l+1/2]^{-1} \delta_{ij}\delta_{lk}.
\end{equation}
The matrix elements given in (B.9) and (B.10) form the orthogonal complete set 
of functions over $M$.The Fourier transform of any square integrable 
function $f\in M$ is given by 
\begin{equation}
f=\sum_{l=0}^\infty \sum_{j=-l}^l [l+1/2] f_j^l D^l_{j0}(G) 
\end{equation}
where the coefficients $f_j^l$  are
\begin{equation}
 f_j^l=\langle f\mid D^l_{j0}(G) \rangle .
\end{equation}

\renewcommand{\theequation}{C.\arabic{equation}}
\setcounter{equation}{0}

{\large C. The Hopf Algebra $U(su_q(2))$}

The Hopf Algebra $U(su_q(2))$ is dual to $A$. Its generated by the elements 
\begin{equation}
{\cal E_\pm}, \ k_\pm=q^{\pm H/4} 
\end{equation}
satisfying  the commutation relations
\begin{equation}
[{\cal E}_{+},{\cal E}_{-}]=\frac{k_{+}^2-k_{-}^2}{q-q^{-1}}, \ \  
 k_{+}k_{-}=k_{-}k_{+}, \ \ \ k_{+}{\cal E}_{+}k_{-}  = q{\cal E}_{+}
\end{equation}
 are the linear functionals on $A$ 
\cite{kn:6}:
\begin {eqnarray}
{\cal E}_1(g)=
\left(
\begin{array}{cc}
0  &  1  \\
0 &  0
\end{array}
\right ),
{\cal E}_{-1}(g)=
\left(
\begin{array}{cc}
0  &  0  \\
1  &  0
\end{array}
\right )  \nonumber \\
k_{\pm}(g)=
\left(
\begin{array}{cc}
q^{\pm 1/2}  &  0  \\
0 &  q^{\mp 1/2}
\end{array}
\right )
\end{eqnarray}  
The extensions of the functionals (C.1) on A are 
\begin {equation}
{\cal E}_\pm (xy)= {\cal E}_\pm (x) k^{+}(y)+k^{-}(x){\cal E}_\pm(y), \ \ \ \  
k^{\pm}(xy)= k^{\pm}(x) k^{\pm} (y)
\end{equation}  
By means of (C.2) and (C.3) we define the differential form of the 
co-representation (B.7) as
\begin {equation}
D(\psi)d^l_{0,j}=\psi(d^l_{k,j})d^l_{0,k}
\end{equation}  
where $\psi$ is one of the elements (C.1). 
We then have
\begin{equation}
D({\cal E_\pm}) d^l_{0j}(M)=([l+1\pm j][l\mp j])^{1/2}d^l_{0,j\pm 1} (M) 
\end{equation}
\begin{equation}
D(k_\pm)  d^l_{0j}(M)=q^{\mp j}d^l_{0j}(M)
\end{equation}
The center of the Hopf algebra $U(su_q(2))$
\begin{equation}
{\cal C} ={\cal E}_{-}{\cal E}_{+} +(\frac{qk_{-}-q^{-1}k_{+}}
{q-q^{-1}})^2
\end{equation}
satisfies the e-value equation
\begin{equation}
(D({\cal C})-[l+1/2]^2)d^l_{0j} =0.
\end{equation}

\vfill
\eject


\begin{thebibliography}{99}
\bibitem{kn:1} \"{O}.F. Dayi and I.H. Duru,
{\em J. Phys. A: Math. Gen.} $\underline{28}$, 2395 (1995); and 
{\em Int. J. Phys.} (in press).
\bibitem{kn:2}N.D. Birrell and P.C.W. Davies,
{\em Quantum Fields in Curved Space,} 
Cambridge Univ. Press, Cambridge, England 1982.
\bibitem{kn:3} See for example I.H. Duru,
{\em Phys.Rev.} $\underline{D30}$, 2121 (1984);
{\em Phys. Lett.} $\underline{119A}$, 163 (1986);
and H. Ahmedov and I.H. Duru,
{\em J. Phys. A: Math. Gen.} $\underline{30,}$ 173 (1997).
\bibitem{kn:4} I.H. Duru and N. \"{U}nal,
{\em Phys. Rev.}$\underline{D34}$, 959 (1986).
\bibitem{kn:5} L.L. Vaksman,Y.S. Soibelman,
{\em Func.Anal.Appl.} $\underline{22}$, 170 (1988);
T. Masuda,K. Mimachi,Y.Nakgami,M. Noumi and K. Ueno,
{\em J. Func. Anal.} $\underline{99}$, 127 (1991);
M. Noumi and K. Mimachi,
{\em in "Lecture Notes in Mathematics"} n.1510, pp.221; P.P. Kulish ed.,
Berlin, Springer 1992;
T.H. Koornwinder,
{\em Proc. Kon. Ned. Akad. Wet., Series A,} $\underline{92}$, 97 (1989);
H.T. Koelink and T.H. Koornwinder,
{\em Proc. Kon. Ned. Akad. Wet., Series A,} 
$\underline{92}$, 443 (1989).
\bibitem{kn:6} N.Ya. Vilenkin and A.O. Klimyk,
{\em Representation  of  Lie Groups and Special Functions}, 
vol.3, Dordrecht;Kluwer Akad. Publ.,1992.
\bibitem{kn:7} F. Bonechi, N. Ciccoli, R. Giachetti, 
E. Sorace and M. Tarlini,
{\em Comm. Math.Phys.}$\underline{175}$, 161 (1996).
\end{thebibliography}
\end{document}